\input amstex 
\documentstyle{amsppt}
\input bull-ppt
\keyedby{bull351/kmt}

\topmatter
\cvol{28}
\cvolyear{1993}
\cmonth{January}
\cyear{1993}
\cvolno{1}
\cpgs{116-122}
\title Best uniform rational approximation of $x^\alpha$ 
on $[0,1]$ \endtitle
\author Herbert Stahl\endauthor
\shorttitle{Uniform rational approximation of $x^\alpha$ 
on $[0,1]$}
\address TFH-Berlin/FB2, Luxemburger Str. 10, D-1000 
Berlin 65, Germany\endaddress
\date April 22, 1992\enddate
\thanks Research supported by the Deutsche 
Forschungsgemeinschaft (AZ: Sta
299/4-2)\endthanks
\abstract A strong error estimate for the uniform rational 
approximation of
$x^\alpha$ on $[0,1]$ is given, and its proof is sketched. 
Let
$E_{nn}(x^\alpha,[0,1])$ denote the minimal approximation 
error in the uniform
norm. Then it is shown that
$$\lim_{n\to\infty}e^{2\pi\sqrt{\alpha
n}}E_{nn}(x^\alpha,[0,1])=4^{1+\alpha}|\sin\pi\alpha|$$
holds true for each $\alpha>0$.\endabstract
\endtopmatter

\document
\heading 1. Introduction\endheading
\par Let $\Pi_n$ denote the set of all polynomials of 
degree at most $n\in\Bbb
N$ with real coefficients; $\Scr R_{mn}$ denote the set 
$\{p/q|p\in\Pi_m$,
$q\in\Pi_n$, $q\nequiv0\}$, $m,n\in\Bbb N$, of rational 
functions; and the {\it
best rational approximant} $r^*_{mn}\in\Scr R_{mn}$, 
$m,n\in\Bbb N$, and the
{\it minimal approximation error} 
$E_{mn}=E_{mn}(x^\alpha,[0,1])$ be defined by
$$E_{mn}\coloneq\|x^\alpha-r^*_{mn}\|_{[0,1]}=\inf_{r\in\Scr
R_{mn}}\|x^\alpha-r\|_{[0,1]},\tag 1.1$$
where $\|\cdot\|_K$ denotes the sup norm on 
$K\subseteq\Bbb R$. It is well
known that the best approximant $r^*_{mn}$ exists and is 
unique within $\Scr
R_{mn}$ (cf. [Me, \S\S9.1, 9.2] or [Ri, \S5.1]). The 
unique existence also
holds in the special case $(n=0)$ of {\it best polynomial 
approximants}.
\par Since $f_\alpha(x)\coloneq|x|^\alpha$ is an even 
function on $[-1,1]$, the
same is true for its unique approximant
$r^*_{mn}=r^*_{mn}(f_\alpha,[-1,1];\cdot)$, and 
consequently a substitution of
$z^2$ by $z$ shows that approximating $|x|^{2\alpha}$ on 
$[-1,1]$ and
$x^\alpha$ on $[0,1]$ poses an identical problem. We have
$$E_{2m,2n}(|x|^{2\alpha},[-1,1])=E_{mn}(x^\alpha,[0,1])%
\quad\text{for all
}m,n\in\Bbb N.\tag 1.2$$
\par From Jackson's and Bernstein's theorems about the 
interdependence of
approximation speed and smoothness of the function to be 
approximated (cf. [Me,
\S\S5.5, 5.6]) we know that in case of $\alpha\in\Bbb R_+
\backslash\Bbb N$ the
minimal error
$E_{m,0}(|x|^\alpha,[-1,1])$ behaves like $O(m^{-\alpha})$ 
as
$m\to\infty$. In [Be1, Be2] Bernstein proved a result that 
is deeper and much
more difficult to obtain; he showed that the limit
$$\lim_{m\to\infty}m^\alpha 
E_{m,0}(|x|^\alpha,[-1,1]):=\beta(\alpha)\tag
1.3$$
exists for each $\alpha>0$; however, an explicit 
expression for the constant
$\beta(\alpha)$, $\alpha>0$, still is not known. 
\par In case $\alpha=1$ the number 
$\beta\coloneq\beta(1)=0.28016\dots$ is
known as {\it Bernstein's constant}. In [Bel, p. 56] 
Bernstein raised the
question whether $\beta$ can be expressed by known 
transcendentals or whether
it defines a new one, and based on numerical upper and 
lower bounds for
$\beta$, which he calculated up to a precision of 
$\pm0.005$, he made the
tentative conjecture $\beta\overset?\to=1/(2\sqrt\pi)$, 
which, however, has
been disproved in [VC1] by extensive and nontrivial high 
precision
calculations.
\par For large values of $\alpha$, Bernstein was able to 
establish in [Be2] an
asymptotic expression. He showed that
$$\lim_{\alpha\to\infty}\frac{\beta(\alpha)}{\Gamma(%
\alpha)|\sin(\pi\alpha/2)|}=\frac1\pi.\tag
1.4$$
\par While Bernstein's investigations on best polynomial 
approximation of $|x|$
and $|x|^\alpha$ were published in the period between 1909 
and 1938, the study
of best rational approximation of $|x|^\alpha$ was only 
started in 1964 by
Newman's surprising (at the time) result in [Ne] that
$$\tfrac12e^{-9\sqrt n}\leq E_{nn}(|x|,[-1,1])\leq 
3e^{-\sqrt n}\quad\text{for
all }n=4,5,\dots.\tag 1.5$$
A comparison of this result with (1.3) shows that the 
convergence behavior of
rational approximants is essentially better than that of 
polynomials.
Newmann's investigation has triggered a whole series of 
contributions, from
which we select a short list with papers that contain 
substantial improvements
of the error estimate in the uniform norm.
$$\alignat2
E_{nn}(x^\alpha,[0,1])\leq&e^{-c(\alpha)\root 3\of 
n},\qquad\alpha\in\Bbb
R_+\qquad&&\text{in [FrSz]};\\
E_{nn}(x^{1/3},[0,1])\leq&e^{-c\sqrt n},\qquad&&\text{in 
[Bu1]};\\
E_{nn}(x^\alpha,[0,1])\leq&e^{-c(\alpha)\sqrt 
n},\qquad\alpha\in\Bbb
R_+,\qquad&&\text{in [Go1]};\\
\tfrac13e^{-\pi\sqrt{2n}}\leq 
E_{nn}(x^{1/2},[0,1])\leq&e^{-\pi\sqrt{2n}(1-\Scr
O(n^{-1/4}))},\qquad&&\text{in [Bu2]};\\
e^{-c(\alpha)\sqrt 
n}\leq&E_{nn}(x^\alpha,[0,1]),\qquad\alpha\in
Q_+\backslash\Bbb N,\qquad&&\text{in [Go2]};\\
e^{-4\pi\sqrt{\alpha n}(1+\varepsilon)}\leq 
E_{nn}(x^\alpha,[0,1])\leq&e^{-\pi
\sqrt{\alpha n}(1-\varepsilon)},\\
\alpha&\in\Bbb R_+\backslash\Bbb N,\ 
\varepsilon>0,n\geq 
n_0(\varepsilon,\alpha)\qquad&&\text{in [Go3]};\\
E_{nn}(x^{1/2},[0,1])\leq&cne^{-\pi\sqrt{2n}},\qquad&&%
\text{in [Vj1]};\\
\tfrac13e^{-\pi\sqrt{2n}}\leq
E_{nn}(x^{1/2},[0,1])\leq&ce^{-\pi\sqrt{2n}}\qquad&&%
\text{in [Vj2]};\\
e^{-c_1(s)\sqrt n}\leq E_{nn}(\root s\of 
x,[0,1])\leq&e^{-c_2(s)\sqrt n},\qquad
s\in\Bbb N\qquad&&\text{in [Tz]}.\endalignat
$$
Here $c$, $c(\alpha),\dots$ denote constants. Relation 
(1.2) allows us to
transfer these results to the problem of approximating 
$|x|^\alpha$ on
$[-1,1]$.
\par It may be appropriate to repeat a remark from [Go2], 
where it was pointed
out that Newman's result can be obtained rather 
immediately from an old result
(from 1877) by Zolotarev. In this sense, the investigation 
of rational
approximants dates back even further than Bernstein's work 
on polynomial
approximants.
\par The best result known so far for the rational 
approximation of $x^\alpha$
on $[0,1]$ was obtained independently by Ganelius [Ga] in 
1979 and by
Vjacheslavov [Vj3] in 1980. They proved that for 
$\alpha\in\Bbb
R_+\backslash\Bbb N$ there exists a constant 
$c_1=c_1(\alpha)>0$ such that
$$\liminf_{n\to\infty}e^{2\pi\sqrt{\alpha 
n}}E_{nn}(x^\alpha,[0,1])\geq
c_1(\alpha)\tag 1.6$$
and, conversely, that for each positive rational number 
$\alpha\in\Bbb Q_+$
there exists a constant $c_2=c_2(\alpha)<\infty$ such that
$$\limsup_{n\to\infty}e^{2\pi\sqrt{\alpha 
n}}E_{nn}(x^\alpha,[0,1])\leq
c_2(\alpha).\tag 1.7$$
Both authors were not able to show that $c_2=c_2(\alpha)$ 
depends continuously
on $\alpha$. Thus, inequality (1.7) remained open for 
$\alpha\in\Bbb
R_+\backslash\Bbb Q$; however, in [Ga] Ganelius was able 
to prove the somewhat
weaker estimate
$$E_{nn}(x^\alpha,[0,1])\leq 
c_2(\alpha)e^{2\pi\sqrt{\alpha n}+c_3(\alpha)\root
4\of n}\quad\text{for }n\geq 
n_0(c_2(\alpha),c_3(\alpha)),\tag 1.8$$
which holds for each $\alpha>0$. In (1.8) $c_2(\alpha)$ 
and $c_3(\alpha)$ are
constants depending on $\alpha$.
\par The results (1.6)--(1.8) give the correct exponent 
$-2\pi\sqrt{\alpha n}$
in the error formula; however, nearly nothing is said 
about the coefficient in
front of the exponential term. The determination of this 
coefficient is the
subject of the present note. Practically as a byproduct, 
we prove the upper
estimate (1.7) for irrational exponents $\alpha\in\Bbb R_+
\backslash\Bbb Q$.
\heading 2. The result\endheading
\thm{Theorem 1} The limit
$$\lim_{n\to\infty}e^{2\pi\sqrt{\alpha
n}}E_{nn}(x^\alpha,[0,1])=4^{1+\alpha}|\sin\pi\alpha|\tag 
2.1$$
holds for each $\alpha>0$.
\ethm
\rem{Remarks} (1) From \kern-.39pt (2.1) \kern-.39pt we 
\kern-.39pt deduce
\kern-.39pt that \kern-.39pt the \kern-.39pt approximation 
\kern-.39pt error
\kern-.39pt $E_{nn}(x^\alpha\!,[0,\!1])$ has the 
asymptotic behavior
$$E_{nn}(x^\alpha,[0,1])=4^{1+
\alpha}|\sin\pi\alpha|e^{-2\pi\sqrt{\alpha n}}
(1+o(1))\quad\text{as }n\to\infty,\tag 2.2$$
and equivalently it follows with (1.2) that
$$E_{nn}(|x|^\alpha,[-1,1])=4^{1+
\alpha/2}|\sin\pi\alpha/2|e^{-\pi\sqrt{\alpha
n}}(1+o(1))\quad\text{as }n\to\infty\tag 2.3$$
for each $\alpha>0$.
\par(2) Not only the explicit expression on the right-hand 
side of (2.1) but
already the existence of the limit represents a result 
difficult to obtain. The
value $4^{1+\alpha}|\sin\pi\alpha|$, $\alpha>0$, is the 
analogue of Bernstein's
constant $\beta(\alpha)$ in (1.3) for the case of rational 
approximation. It
has already been noted in the introduction that an 
explicit expression for
$\beta(\alpha)$ is still not known. The best we know is 
Bernstein's asymptotic
formula (1.4). Since in (1.4) we have considered 
approximation on $[-1,1]$, the
counterpart of the asymptotic value 
$\tfrac1\pi\Gamma(\alpha)|\sin\pi\alpha/2|$
for $\beta(\alpha)$ is the value $4^{1+
\alpha/2}|\sin(\pi\alpha/2)|$ in case of
rational approximation.
\par(3) If we turn our attention to the special case of 
rational approximation
of $|x|$ on $[-1,1]$, then it follows from (1.2) that
$$E_{2n,2n}(|x|,[-1,1])=E_{nn}(\sqrt 
x,[0,1])\quad\text{for }n\in\Bbb N,\tag
2.4$$
and hence we deduce from (2.1) that
$$\lim_{n\to\infty}e^{-\pi\sqrt n}E_{nn}(|x|,[-1,1])=4^{1+
1/2}|\sin\tfrac\pi
2|=8.\tag 2.5$$
Limit (2.5) has been conjectured in [VRC] on the basis of 
high precision
calculations and was proved in [St].
\par It may be surprising that in case of rational 
approximation, which is in
many respects more complex than the polynomial case, limit 
(2.5) has a rational
value, while in the polynomial case Bernstein's question 
in [Be1] about the
character of the number $\beta=\beta(1)$ is still open and 
the numerical
results in [VC1] show that $\beta$ cannot be a rational 
number with a
moderately small denominator.
\endrem
\heading 3. Outline of the proof of Theorem 1\endheading
\par If limit (2.1) is proved for one of the paradiagonal 
sequences
$$\{E_{n+k,n}(x^\alpha,[0,1])\}^\infty_{n=|k|},$$
$k\in\Bbb Z$ fixed, then it
holds also for the diagonal sequence 
$$\{E_{nn}(x^\alpha,[0,1])\}^\infty_{n=1}.$$
It turns out that
$$n+k=m\coloneq n+1+[\alpha],\qquad\alpha\in\Bbb R_+
\backslash\Bbb N,n\in\Bbb
N,\tag 3.1$$
is a good choice for the numerator degree $m$. From the 
theory of best rational
approximants we learn that the error function
$$e_n(z)\coloneq z^\alpha-r^*_{mn}(z),\qquad z\in\Bbb 
C\backslash\Bbb R_-,\tag
3.2$$
has exactly $2n+2+[\alpha]$ zeros in the interval $(0,1)$. 
Hence the theory of
multipoint Pad\'e approximants is applicable, and it gives 
us rather precise
information about the structure of the numerators and 
denominators of
$r^*_{mn}$ (cf. [GoLa; StTo, \S\S6.1, 6.2]).
\par In the next step the error function $e_n$ and the 
approximant $r^*_{mn}$
will be transformed in such a way that the resulting 
function $\Psi_n$ is
analytic in $\Bbb C\backslash\Bbb R$ with possible 
exceptions in a disc
$\Delta(R)$ with radius $R>0$ around the origin. The 
function $\Psi_n$ has
boundary values from both sides of $\Bbb 
R\backslash\overline{\Delta(R)}$ that
allow a comparison with a special logarithmic potential. 
The potential will be
introduced after Theorem 2 below.
\par The transformation of $e_n$ and $r^*_{mn}$ into the 
function $\Psi_n$ is
carried out in several steps. The intermediate functions 
are defined as
$$\gather
r_n(z)\coloneq\frac{z^\alpha-r^*_{mn}(z)}{z^\alpha+
r^*_{mn}(z)};\tag 3.3\\
R_n(w)\coloneq\frac{4w^{2\alpha}-1}{w^\alpha}r_n(%
\varepsilon^{1/\alpha}_nw)-\frac
1{w^\alpha},\qquad\varepsilon_n\coloneq 
E_{mn}(x^\alpha,[0,1]);\tag
3.4\\
\Phi_n(w)\coloneq\frac1{8w^\alpha}\left(R_n(w)+
\sqrt{R_n(w)^2-4}\right);\tag 3.5\\
\Psi_n(w)\coloneq\cases\psi(\Phi_n(w))&\text{for 
}\operatorname{Im}(w)\geq 0,\\
\overline\psi(\Phi_n(w))&\text{for 
}\operatorname{Im}(w)<0,\endcases\tag 3.6
\endgather
$$
with
$$\psi(z)\coloneq\frac 
z{\sin\pi\alpha-i(\cos\pi\alpha)z}.\tag 3.7$$
In (3.4) a new variable $w$ is introduced implicitly by
$$w\coloneq\varepsilon^{-1/\alpha}_nz,\qquad z\in\Bbb 
C.\tag 3.8$$
\par The properties of each new function $r_n,R_n,\Phi_n$ 
have to be studied
carefully. The properties of the last function $\Psi_n$ is 
summerized in
\thm{Lemma 1} The function $\Psi_n$ is analytic and 
different from zero in
$\Bbb C\backslash(\Bbb R\cup\overline{\Delta(R)})$, $R>0$ 
appropriately chosen,
and there exist constants $c_1,\dots,c_4$ such that
$$\gather
|\log|\Psi_n(w)\|\leq c_1|w|^{-2\alpha}\quad\text{for 
}w\in\Bbb
R_-\backslash\Delta(R),n\geq n_0(c_1,R);\tag 3.9\\
|\log|\Psi_n(w)w^\alpha 
4\sin\pi\alpha\|\leq c_2|w|^{-\alpha}
\quad\text{for}\ w\in\Bbb R_+\backslash\Delta(R),n\geq 
n_0(c_2,R);\tag 3.10\\
|\log|\Psi_n(w)\|\leq c_3\quad\text{for 
}w\in\partial\Delta(R),n\geq
n_0(c_3,R).\tag 3.11\endgather$$
If we consider the representation
$$\log|\Psi_n(w)|=\psi_n(w)-\int g_D(w,t)\,d\mu_n(t),\tag 
3.12$$
where $D\coloneq\Bbb C\backslash(\Bbb 
R_-\cup\overline{\Delta(R)})$, $g_D(w,t)$
the Green function in $D$, $\psi_n$ a harmonic function in 
$D$ with
$$\psi_n(w)=\log|\Psi(w)|\quad\text{for }w\in\partial 
D,\tag 3.13$$
and $\mu_n$ a measure with
$$\supp(\mu_n)\subseteq[R,\infty]\qquad\text{and}\qquad%
\mu_n\geq 0\quad\text{on
}[R,\varepsilon^{-1/\alpha}_n],\tag 3.14$$
then in addition to $(3.9)\text{--}(3.11)$ we have the 
inequalities
$$\gather
\|\mu_n|_{[\varepsilon^{-1/\alpha},\infty]}\|\leq 1,\tag 
3.15\\
|\mu_n([R,\varepsilon^{-1/\alpha}_n])-2n|\leq 
c_4\quad\text{for }n\geq
n_0(c_4,R).\tag 3.16\endgather
$$
\ethm
\rem{Remark} Estimates (3.9), (3.10), and (3.16) contain 
the information that
is most relevant for the proof of Theorem 1.
\endrem
\par The function $\log|\Psi_n|$ will be compared with a 
special Green
potential $p_n$. The definition of this potential is based 
on the following.
\thm{Theorem 2} For each $a\in(1,\infty)$ there uniquely 
exists a Green
potential
$$p_a(w)\coloneq\int g_D(w,t)\,d\nu_a(t),\qquad 
D\coloneq\Bbb C\backslash\Bbb
R_-,\tag 3.17$$
with
$$\nu_a\geq 
0,\qquad\operatorname{supp}(\nu_a)=[b_a,a],\qquad 
0<b_a<a,\tag
3.18$$
that satisfies
$$p_a(w)\cases=\log w&\text{for }w\in[b_a,a],\\
>\log w&\text{for }w\in(0,b_a).\endcases\tag 3.19$$
For the constant $b_a$ and the measure $\nu_a$ appearing in
$(3.17)\text{--}(3.19)$ we have
$$b_a\to\sqrt 2\tag 3.20$$
and
$$\tfrac1ae^{\pi\sqrt{2\|\nu_a\|}}\to 4\quad\text{as 
}a\to\infty.\tag 3.21$$
\ethm
\par A proof of Theorem 2 can be derived from [St, Theorem 
2]. It is necessary
to change the domain of definition $D$ by the 
transformation $w\mapsto 1/\sqrt
w$. Basic tools in the proof of Theorem 2 are estimates 
for certain elliptical
integrals.
\par The potential $p_n$ for a comparison with the 
function $\log|\Psi_n|$ is
now defined as
$$p_n(w)\coloneq-\alpha p_a(cw)\tag 3.22$$
with
$$a\coloneq\left|\frac
4{\varepsilon_n}\sin\pi\alpha\right|^{1/\alpha}\quad%
\text{and}\quad
c\coloneq|4\sin\pi\alpha|^{1/\alpha}.\tag 3.23$$
From (3.21), (3.23), and (3.16) together with other 
information provided by
Lemma 1 it then follows that
$$\left[\left|\frac{\varepsilon_n}{4\sin\pi\alpha}%
\right|^{1/\alpha}e^{\pi\sqrt{2
(2n)/\alpha}}\right]^\alpha=\frac{\varepsilon_n}{4|\sin\pi%
\alpha|}
e^{2\pi\sqrt{\alpha n}}\to 4^\alpha\quad\text{as 
}n\to\infty.\tag 3.24$$
Theorem 1 then immediately follows from (3.24).
\par What we have given here is only a sketch of the 
overall structure of the
proof of Theorem 1. Some of the steps demand rather subtle 
analysis.
\rem{Note Added in Proof} In [VC2] Varga and Carpenter 
have calculated
numerical approximations for the right-hand side of (2.1) 
for the six values
$\alpha=j/8$, $j=1,2,3,5,6,7$, and conjectured formula 
(2.1) on the basis of
these numerical values. The numerical results for the two 
cases $\alpha=1/4$
and $\alpha=3/4$ are especially interesting since here the 
right-hand side of
(2.1) is rational.
\endrem

\enddocument